\documentclass[12pt]{article}
\usepackage{latexsym, amsmath, amssymb, amsfonts, amscd}
\usepackage{graphics}
\usepackage[dvips]{epsfig}
\newtheorem{thm}{Theorem}[section]
\newtheorem{prop}[thm]{Proposition}
\newtheorem{lemma}[thm]{Lemma}
\newtheorem{cor}[thm]{Corollary}
\newtheorem{dfn}[thm]{Definition}

\newcommand{\oper}[2]{\newcommand{#1}{\mathop{\mathrm{#2}}\nolimits} }
\oper{\DGA}{DGA} \oper{\END}{End} \oper{\nij}{Nij} \oper{\jac}{Jac}

\newcommand{\call}{{\cal L}}

\def\qed{\rule{2.3mm}{2.3mm}}

\newcommand{\half}{\frac{1}{2}}

\def\lcf{\lbrack\! \lbrack}
\def\rcf{\rbrack\! \rbrack}

\newcommand{\CC}{\mathbb{C}}

\newcommand{\RR}{\mathbb{R}}

\def\cl{{\cal L}}

\def\om{\omega}

\def\oL{\overline L}
\def\oE{\overline E}

\def\oE{\overline E}

\newcommand{\lie}[1]{\mathfrak{#1}}
\newcommand{\vspan}[1]{\langle {#1} \rangle}

\def\lf{{\cal{L}}_{{F}}}

\newcommand{\bra}[2]{\{ #1 \cdot #2\}}

\newcommand{\lbra}[2]{\lcf #1, #2 \rcf}

\newcommand{\inner}[2]{\langle #1, #2 \rangle}

\newcommand{\bproof}{\noindent{\it Proof: }}
\newcommand{\eproof}{\hfill \qed. \vspace{0.2in}}

\newcommand{\Lie}[1]{\operatorname{\textsl{#1}}}

\newcommand{\GL}{\Lie{GL}}

\newcommand{\SO}{\Lie{SO}}

\newcommand{\SP}{\Lie{Sp}}

\newcommand{\SU}{\Lie{SU}}
\newcommand{\su}{\lie{su}}

\begin{document}
\title{Generalized Contact Structures}
\author{ Y.~S.  Poon \thanks{ Address:
    Department of Mathematics, University of California at Riverside,
    Riverside CA 92521, U.S.A., Email: ypoon@ucr.edu. Partially supported by UC-MEXUS
    and NSF-0906264} \
 and   A\"issa Wade \thanks{Address: Department of Mathematics, The Pennsylvania State University,
University Park, PA 16802, U.S.A., Email: wade@math.psu.edu}  }

%\classno{53D (primary), 17B, 32G, 58H (secondary)}

\maketitle
\begin{abstract}  We study  integrability of generalized almost contact structures,
and find conditions under which the main associated  maximal
isotropic vector bundles form Lie bialgebroids. These conditions
differentiate the concept of generalized contact structures from a
counterpart of generalized complex structures on odd-dimensional
manifolds. We name the latter strong generalized contact structures.
Using a Boothby-Wang construction bridging symplectic structures and
contact structures, we find examples to demonstrate that, within the
category of generalized contact structures, classical contact
structures have non-trivial deformations. Using deformation theory
of Lie bialgebroids, we construct new families of strong generalized
contact structures on the three-dimensional Heisenberg group and its
co-compact quotients.
\end{abstract}

\section{Introduction}

The theory of generalized complex
structures is a geometric framework unifying both complex structures and symplectic structures
  \cite{Marco}  \cite{Hitchin}. It is applicable only to even-dimensional manifolds.
  A key feature of this theory is to allow deformation between
  complex and symplectic structures. There are indeed non-trivial examples of such phenomenon on compact
  manifolds \cite{Marco} \cite{Poon}. This phenomenon is a
departure from Moser's theorem on the rigidity of symplectic structures with respect to diffeomorphisms \cite{Moser}.

 For decades, symplectic structures and contact structures have often been
  studied in parallel, beginning as frameworks for
 classical mechanics \cite{Arnold} \cite{Lie}. For instance, both symplectic and contact structures have ``standard" local models.
 Moser's theorem has its counterpart for contact structures \cite{Gray}.
 Boothby and Wang showed that when a contact structure is represented by
 a ``regular" one-form, the underlying manifold is
 foliated and the leave space has a symplectic structure.
 Conversely, the total space  of  a $\SO(2)$-bundle
 on a symplectic manifold with the given symplectic form as a curvature
 form has a contact structure \cite{BW}. From the viewpoint of
 G-structures, contact structures are also related to complex structures.
 This perspective leads to Sasaki's introduction of normal almost
 contact structures on odd-dimensional manifolds \cite{Sasaki1}.

 While much of the similarity between symplectic and contact structures are emphasized,
 and relation between complex
 and contact structures are developed,
 we often ignore a fundamental distinction of contact structures.
 Namely from a  G-structure perspective, symplectic structures and complex structures are integrable.
 The pseudogroup of their
 local models are transitive \cite{Chern} \cite{Chern66}.
 Although the pseudogroup of contact transformations remains to be transitive \cite{BG},  contact structures
 are not integrable G-structures.

 In this article, we continue a recent quest for developing an analogue of generalized complex structures on
 odd-dimensional manifolds \cite{IW} \cite{PW} \cite{Va06}.  In \cite{IW}, the second author and her collaborator
developed a concept called generalized almost contact structures.
 Its development is based on the theory of Dirac
structures and 1-jet bundles of the underlying manifolds. In
\cite{Va06}, Vaisman developed the concept of  ``generalized almost
contact structures of co-dimension k". He also introduced  and studied generalized F-structures and
CRF-structures in \cite{Va08}. We focus on co-dimension one case in Vaisman's development,
and simply call it a ``generalized almost contact structure" in this article (see Definition \ref{restricted}).

In a recent paper, we investigate the integrability of such
structures from Sasaki's perspective \cite{PW}. However, the theory
of generalized complex structures is developed in the context of Lie
bialgebroids \cite{Marco} \cite{LWX}. This concept requires the
splitting of the complexification of the direct sum of the tangent
bundle $TN$ and the cotangent bundle $T^*N$ over an even-dimensional
manifold $N$ into the direct sum of two maximally isotropic
subbundles $L$ and $L^*$. Integrability is in terms of the
closedness of the spaces of sections of these bundles under the
Courant bracket  \cite{Marco}. The pair $L$ and $L^*$ forms a Lie
bialgebroid.  The core of this paper is to analyze generalized
almost contact structures in such context.

Given a generalized almost contact structure ${\cal J}$ on a
manifold $M$, readers will see that the bundle $(TM\oplus
T^*M)_\mathbb{C}$ splits into the direct sum of two maximally
isotropic sub-bundles, $L$ and its dual $L^*$.   Unlike generalized
almost complex structures on even-dimensional manifolds, $L^*$ is
not complex-conjugate linearly isomorphic to $L$. Therefore, one has
to analyze $L$ and $L^*$ individually. In Section \ref{obstruction
section},  we identify the obstruction for $\Gamma(L^*)$ to be
closed under the assumption that the space $\Gamma(L)$ of sections
of $L$ is closed under the Courant bracket. The main result of this
section is Theorem \ref{bi}.

If the structure ${\cal J}$ is defined by a classical contact 1-form,
we show that  the space $\Gamma(L)$
 is closed under the Courant bracket (see Section \ref{classical contact}).
However, by analyzing the local model of a contact structure, we
find that  the obstruction for $\Gamma(L^*)$ to be closed under the
Courant bracket does not vanish (see Proposition \ref{contact}).
Therefore, we define a generalized (integrable) contact structure to
be a generalized almost contact structures whose corresponding
$\Gamma(L)$ is closed under the Courant bracket, while $\Gamma(L^*)$
is not necessarily closed (see Definition \ref{closed}). When both
$\Gamma(L)$ and $\Gamma(L^*)$ are closed, we consider the given
structure ${\cal J}$ as a natural counterpart of a generalized
complex structure on an odd-dimensional manifold, and name ${\cal
J}$ a strong generalized contact structure  (see Definition
\ref{odd}). The distinction between generic generalized contact
structures and the strong version could be conceived as an extension
of the fact that classical contact structures are not integrable
G-structures.

We find examples for these new concepts from two different sources.
One is from classical geometry. Another
is through deformation theory.

The classical analogues of symplectic and complex structures were
discovered nearly half a century ago. When studying infinitesimal
automorphisms of symplectic structures, Libermann developed the
concept of almost cosymplectic structures \cite{Libermann}. As a
G-structure, it is a reduction of the structure group of a
$(2n+1)$-dimensional manifold from $\GL(2n+1, \RR)$ to
$\{1\}\times\SP(n, \RR)$ \cite{Fuji} \cite{Libermann}. In terms of
tensors, it is equivalent to  the choice of a 1-form $\eta$ and a
2-from $\theta$ such that $\eta\wedge\theta^n\neq 0$ at every point
of the manifold.  An almost cosymplectic structure $(\eta, \theta)$
is a cosymplectic structure if it is an integrable G-structure. It
is equivalent to both $\eta$ and $\theta$ being closed forms
\cite{Libermann}.  By choosing $\theta=d\eta$, it is immediate that
 a contact 1-form $\eta$ yields an
almost cosymplectic structure, but it  is \em not \rm integrable
as  $d\eta$ is non-zero everywhere.

Treating an almost complex structure on a $2n$-dimensional manifold
$N$ as a reduction of the principal bundle of frames from $\GL(2n,
\RR)$ to $\GL(n, \CC)$, we obtain an $(1,1)$-tensor $J$ on the
manifold $N$ such that $J\circ J=-\mathbb{I}$. An almost contact
structure on an odd-dimensional manifold $M$ is a triple
 $(F, \eta, \varphi)$ consisting of a vector field $F$, a one-form $\eta$ and a $(1,1)$-tensor
 $\varphi$ such that $\varphi^2=-\mathbb{I}+F\otimes \eta$.
This triple could be used to define naturally an almost complex structure on the cone over $M$
 \cite{Blair} \cite{BG} \cite{Sasaki1}.
 When this almost complex structure is integrable,
 the triple is called a normal almost contact structure.
 Readers are
 warned of the very unfortunate historical fact that without an
 auxiliary geometric
 object, a contact 1-form does not
 naturally define an almost contact structure, nor does a normal
 almost
 contact structure \cite{BG}.

In Section \ref{almost cosym} and Section \ref{almost contact}, it
is respectively shown that cosymplectic structures and normal almost
contact structures are examples of strong generalized contact
structures. Since such structures are associated with Lie
bialgebroid theory, we are able to apply a deformation theory as
developed in \cite{LWX} to generate new and non-classical examples
(see Section \ref{H3}).

As noted in  Section \ref{classical contact}, classical contact
structures are examples of non-strong generalized contact
structures. We illustrate this point on $\SU(2)$ in Section
\ref{SU2}. To find non-trivial new examples, we apply a Boothby-Wang
construction on a $\SO(2)$-bundle on the Kodaira surface $N$
\cite{BW}. We first note that there exists an analytic family of
generalized complex structures $J_t$ with parameter $t$ such that
$J_0$ is a classical complex structure on a Kodaira surface $N$, and
$J_1$ is a symplectic structure on $N$ \cite{Poon}. Following
\cite{BW}, we construct a family of generalized contact structure
${\cal J}_t$ on a principal $\SO(2)$-bundle $M$ over the Kodaira
surface such that ${\cal J}_1$ is associated to a classical contact
1-form on $M$. The quotient of the structure ${\cal J}_t$ on $M$  by
the fundamental vector field of the principal bundle yields the
family $J_t$ on $N$. This example explicitly illustrates that
classical contact 1-forms, when conceived as generalized contact
structures, are not necessarily rigid. It is a departure from Gray's
theorem \cite{Gray}.

\vspace{0.25in}

\noindent{\bf Acknowledgement.} We thank Charles Boyer, Camille
Laurent-Gengoux, Jean-Pierre Marco,  and Pol Vanhaecke for useful
discussions. The first author thanks the hospitality of Centro de
Investigacion en Matematicas (C.I.M.A.T.) in Guanajuato, Mexico, and
the financial support of UC-MEXUS and NSF DMS-0906264. We also thank
Izu Vaisman, the referee and Andrew Swann of LMS for very helpful suggestions.

\section{Generalized Contact Structures}

For a manifold $M$ of any dimension, consider the vector bundle
$TM\oplus T^*M \to M$.
 Its space of sections is endowed with
two natural $\mathbb R$-bilinear operations.
\begin{itemize}
 \item A symmetric bilinear form $\langle \cdot, \cdot\rangle$ is
 defined by
\begin{equation}\label{symmetric pairing}
 \langle X+ \alpha,  Y+\beta\rangle=\frac12( \iota_X\beta+\iota_Y\alpha).
 \end{equation}
\item  The Courant bracket is given by
\begin{equation}\label{courant bracket}  \lcf X+ \alpha, Y+\beta \rcf= [X, Y]+{\call}_X\beta- {\call}_Y\alpha-
 \frac{1}{2}d(\iota_X\beta- \iota_Y\alpha).
 \end{equation}
\end{itemize}

We adopt the notations: $(\pi^\sharp\alpha)(\beta)=\pi(\alpha,
\beta)$, and $Y(\theta^\flat X)$$=\theta(X,Y)$ for any 1-forms
$\alpha$ and $ \beta$, 2-form $\theta$, bivector field $\pi$, and
vector fields $X$ and $Y$.

The bundle $TM\oplus T^*M $ with the non-degenerate pairing $\langle
-, -\rangle$ in (\ref{symmetric pairing}) and Courant bracket
(\ref{courant bracket}) above form a fundamental example of {\em
Courant algebroid}  \cite{Courant} \cite{LWX}. The natural
projection $\rho$ from the direct sum to the summand $TM$ is called
the anchor map.

We will consider the complexified bundles, and complex-linearly
extend the symmetric bilinear form and the Courant bracket to obtain
complex Courant algebroids.

\subsection{Generalized  almost contact structures}\label{general}

\begin{dfn}{\rm \cite{PW} \cite{Va06}} \label{restricted} A generalized  almost contact pair on
a smooth odd-dimensional manifold $M$
consists of a bundle endomorphism $\Phi$ from $TM \oplus T^*M$ to
itself and a section ${F}+\eta$ of $TM \oplus T^*M$ such that  $\Phi+\Phi^*=0$,
$\eta({F})=1$, $\Phi(F)=0$, $\Phi(\eta)=0$, and $\Phi\circ\Phi=-\mathbb{I}+{F}\odot\eta.$
\end{dfn}
Here $F\odot\eta$ acts on $TM\oplus
T^*M$ as a symmetric bundle endomorphism, i.e. when $X+\alpha$ is a section of
$TM\oplus T^*M$, then as a matter of definition,
\[
(F\odot\eta)(X+\alpha):=\eta(X) F+\alpha(F) \eta.
\]

The pair of tensors $(\Phi, F+ \eta)$ is equivalent to another
 pair $(\Phi', F'+ \eta')$ if there exists a function $f$ without
zero on the manifold $M$ such that
\begin{equation}
\Phi'=\Phi, \quad \eta'=f\eta, \quad F'=\frac{1}{f}F.
\end{equation}

\begin{dfn} \label{restricted 2} A generalized almost contact
structure on $M$ is an equivalent class of such pair $(\Phi, F + \eta)$.
\end{dfn}

In terms of components, a generalized almost contact structure is
given by an equivalent class of tensorial objects: ${\cal J}=(F,
\eta, \pi, \theta, \varphi)$ where ${F}$ is a vector field, $\eta$ a
1-form, $\pi$ a bivector field, $\theta$ a 2-form, and $\varphi$ a
(1,1)-tensor. They are subjected to the following relations.
\begin{eqnarray}
&\theta^\flat\varphi=\varphi^*\theta^\flat, \quad
 \varphi\pi^\sharp=\pi^\sharp\varphi^*,& \label{Relation1}\\
&\varphi^2+\pi^\sharp\theta^\flat=-\mathbb{I}+{F}\otimes \eta,
\mbox{ and }
(\varphi^*)^2+\theta^\flat\pi^\sharp=-\mathbb{I}+\eta\otimes {F}.&
\label{Relation2}\\
 &\eta\circ \varphi=\varphi^*\eta=0,  \eta \circ\pi^\sharp=\pi^\sharp\eta=0,  \iota_{F}\varphi=0,
  \iota_{F}\theta=0, \iota_{F}\eta=1.& \label{Relation3}
\end{eqnarray}
The bundle map $\Phi: TM \oplus T^*M \rightarrow TM \oplus T^*M$ is
given by
\begin{equation}\label{Phi in general}
\Phi= \left(
\begin{array}{cc}
\varphi  & \pi^\sharp \\
\theta^\flat & - \varphi^*
\end{array}
\right).
\end{equation}

\subsection{The associated complex vector subbundles}\label{bundles}
 Consider the above  bundle map $\Phi$. It has one real eigenvalue, namely $0$. The
corresponding eigenbundle is trivialized by $F$ and $\eta$
respectively. We denote these bundles by $L_F$ and $L_\eta$.
 Let $\ker \eta$ be the distribution on the
manifold $M$ defined by the point-wise kernel of the 1-form $\eta$.
Similarly, $\ker F$ is the subbundle of $T^*M$ defined by the
point-wise kernel of the vector field $F$ with respect to its
evaluation on differential 1-forms. On the complexified bundle
$(TM\oplus T^*M)_\CC$, $\Phi$ has three eigenvalues, namely ``$0$",
``$+i$" and ``$-i$".
  Define
\begin{eqnarray*}
 &E^{(1,0)}=\{ e - i \ \Phi (e) \ | \ e \in \ker
\eta \oplus \ker {F} \},&\\
 &\hspace{1in} E^{(0,1)}=\{ e + i\ \Phi (e) \ | \ e \in \ker
\eta \oplus \ker {F} \}.&
\end{eqnarray*}
Then $L_F\oplus L_\eta$ is the $0$-eigenbundle,  $ E^{(1,0)}$ is the
$+i$-eigenbundle and  $E^{(0,1)}$ is the $-i$-eigenbundle. We have a
natural splitting: $ (TM \oplus T^*M)_\CC = L_{{F}} \oplus L_{\eta}
\oplus E^{(1,0)} \oplus E^{(0,1)}. $ It is apparent that this
decomposition does not depend on any choice of representatives
within an equivalent class of generalized contact forms. A choice of
a trivialization of a real subbundle $L_F$ in $TM$,  a
trivialization of its dual $L_\eta$ in $T^*M$, and the subsequent
choice of $E^{(1,0)}\oplus E^{(0,1)}$ determines a generalized
contact pair.

In subsequent analysis, the following four different complex vector
bundles will play different roles. Namely,
\begin{eqnarray}
&L:=L_{F} \oplus E^{(1,0)} , \quad {\oL}:=L_{F}\oplus E^{(0,1)},& \nonumber\\
  &\hspace{1in} L^*:=L_\eta \oplus E^{(0,1)}, \quad {\oL}^*:=L_\eta\oplus E^{(1,0)}.&\label{fundamental}
\end{eqnarray}
As $L_{F}$ is the complexification of a real line bundle, its
conjugation is itself. Therefore, the complex conjugation map sends
$L$ to $\oL$. On the other hand, through the symmetric pairing
(\ref{symmetric pairing}), $L^*$ is complex-linearly isomorphic to
the dual of $L$. All these bundles are independent of choice of representatives
of a generalized almost contact structure.

\begin{lemma}\label{isotropic} The bundles $E^{(1,0)}$, $E^{(0,1)}$, $L$,
$\overline{L}$, $L^*$ and ${\overline{L}}^*$ are isotropic with respect to
 the symmetric pairing $\langle -,
-\rangle$.
\end{lemma}
\bproof Suppose that $X+\alpha$ is section of
$\ker\eta\oplus\ker{F}$. Then
\[
\Phi(X+\alpha)=\varphi(X)+\pi^\sharp(\alpha)+\theta^\flat(X)-\varphi^*(\alpha).
\]
By constraints (\ref{Relation3}), $\Phi(X+\alpha)$ is again a
section of $\ker\eta\oplus\ker{F}$. Therefore,
\begin{equation}\label{f and eta}
\langle {F}, \Phi(X+\alpha)\rangle=0 \quad \mbox{ and } \quad
\langle \eta,\Phi(X+\alpha)\rangle=0.
\end{equation}

If both $X+\alpha$ and $Y+\beta$ are sections of
$\ker\eta\oplus\ker{F}$, then
\begin{eqnarray*}
&&\langle X+\alpha-i\Phi(X+\alpha), Y+\beta-i\Phi(Y+\beta)\rangle\\
&=&\langle X,\beta\rangle-i\langle X,
\theta^\flat(Y)-\varphi^*(\beta)\rangle  -i\langle \alpha,
\varphi(Y)+\pi^\sharp(\beta)\rangle\\
&+&\langle Y, \alpha\rangle-i\langle Y,
\theta^\flat(X)-\varphi^*(\alpha)\rangle  -i\langle \beta,
\varphi(X)+\pi^\sharp(\alpha)\rangle\\
&-&\langle \varphi(X)+\pi^\sharp(\alpha),
\theta^\flat(Y)-\varphi^*(\beta) \rangle -\langle
\varphi(Y)+\pi^\sharp(\beta), \theta^\flat(X)-\varphi^*(\alpha)
\rangle.
\end{eqnarray*}
Since $\theta$ and $\pi$ are skew-symmetric, the above is reduced to
\begin{eqnarray*}
&=&\langle X,\beta\rangle+\langle \alpha, Y\rangle\\
&-&\langle \varphi(X)+\pi^\sharp(\alpha),
\theta^\flat(Y)-\varphi^*(\beta) \rangle -\langle
\varphi(Y)+\pi^\sharp(\beta), \theta^\flat(X)-\varphi^*(\alpha)
\rangle.
\end{eqnarray*}
By constraints (\ref{Relation2}), it is further reduced to
\[
-\langle \varphi(X), \theta^\flat(Y) \rangle +\langle
\pi^\sharp(\alpha), \varphi^*(\beta) \rangle
 -\langle \varphi(Y),
\theta^\flat(X) \rangle
 +\langle \pi^\sharp(\beta),
\varphi^*(\alpha) \rangle.
\]
By (\ref{Relation1}), it is equal to zero. It follows that $E^{(1,0)}$ is isotropic.

Taking complex conjugation, we find that $E^{(0,1)}$ is also
isotropic. By (\ref{f and eta}),  the pairings between sections of
$L_F$ or $L_\eta$ with those of $E^{(1,0)}\oplus E^{(0,1)}$ are
always equal to zero. Therefore $L=L_{F}\oplus E^{(1,0)}$ is
isotropic. A similar computation shows that $L^*=L_\eta\oplus
E^{(0,1)}$ is isotropic. \eproof

\begin{dfn}\label{closed}
Given a generalized almost contact structure, if the space $ \Gamma(L)$ of sections of  the associated bundle $L$ is closed under the Courant bracket  then
 the generalized almost contact structure is  simply called a generalized contact structure.
\end{dfn}

 Since $L_{F}$ is a rank-1 bundle, it is apparent that $\lcf
\Gamma(L_{F}),  \Gamma(L_{F})\rcf\subseteq  \Gamma(L_{F})$. Therefore, the non-trivial
conditions  for $\lcf \Gamma(L), \Gamma(L) \rcf\subseteq \Gamma(L)$ are due to the following
two inclusions.
\begin{eqnarray*}
&\lcf \Gamma\left(L_{F}\right), \Gamma\left(E^{(1,0)}\right)\rcf
\subseteq \Gamma\left(L_{F} \oplus
E^{(1,0)}\right),&\\
 &\hspace{1in}\lcf \Gamma\left(E^{(1,0)}\right),
\Gamma\left(E^{(1,0)}\right)\rcf \subseteq  \Gamma\left(L_{F} \oplus
E^{(1,0)}\right).&
\end{eqnarray*}

As a consequence of Lemma \ref{isotropic},  all four bundles given
in (\ref{fundamental}) are {\em maximally} isotropic with respect to
the pairing (\ref{symmetric pairing}) in $(TM\oplus T^*M)_\CC$.
Combined with the concept given in Definition \ref{closed}, the
definition of ``Dirac structures" \cite{Courant}, and the definition
of ``quasi"-Lie bialgebroid in \cite{Roy}, we have

\begin{cor} When ${\cal J}=(F, \eta, \pi, \theta,
\varphi)$ represents a generalized contact structure, the associated
bundle $L$ is a Dirac structure. In addition, the bundle $L^*$ is a
transversal isotropic complement of $L$ in the Courant algebroid
$\left((TM\oplus T^*M)_\CC, \langle -, -\rangle, \lcf -,
-\rcf\right)$. In other words, the pair $L$ and  $L^*$ is a
quasi-Lie bialgebroid.
\end{cor}

\subsection{Obstruction to integrability of the dual bundle $L^*$}\label{obstruction section}

A lack of natural isomorphism between $L$ and $L^*$ means that when
$\Gamma(L)$ is closed under the Courant bracket, $\Gamma(L^*)$ is not
necessarily closed. It is a major departure from the theory of
generalized complex structures on even-dimensional manifolds. In
this section, with \cite{LWX} and \cite{Marco} as our key
references, we examine the obstruction for both $L$ and $L^*$ being
closed.

Recall \cite{Mac-2} that a complex Lie algebroid on a manifold $M$
is a complex vector bundle $V$ together with a bundle map $\rho:
V\to TM_\CC$, called the anchor map, and a bracket $\lbra{-}{-}$ on
the space of sections of $V$ such that for any sections $s_1, s_2,
s_3$ of $V$, and any smooth function $f$ on $M$,

\begin{itemize}
\item $\lbra{s_1}{s_2}=-\lbra{s_2}{s_1}$;
\item
$\lbra{\lbra{s_1}{s_2}}{s_3}+\lbra{\lbra{s_2}{s_3}}{s_1}+\lbra{\lbra{s_3}{s_1}}{s_2}=0$;
\item $\lbra{s_1}{fs_2}=f\lbra{s_1}{s_2}+\big(\rho(s_1)f\big)s_2$;
\item  $\rho(\lbra{s_1}{s_2})=[\rho(s_1), \rho(s_2)]$.
\end{itemize}

In our previous discussion on generalized contact structures, we
focus on the bundle $L$. In terms of Lie algebroid, the inclusion of
$L$ in $(TM\oplus T^*M)_\CC$, followed by the natural projection
onto the first summand is an anchor map. When $\Gamma(L)$ is closed
under the Courant bracket, the restriction of the Courant bracket to
$L$ completes a construction of a Lie algebroid structure on $L$.

Given the natural anchor map $\rho$ on $(TM\oplus T^*M)_\CC$ and
assume that $\Gamma(L)$ is closed, the next issue is whether the
space $\Gamma(L^*)$ of sections of $L^*$ is closed under the Courant
bracket.

It is known that the obstruction for $\Gamma(L^*)$ to be closed is
due to an alternating form on $L^*$ \cite[Lemma 3.2]{LWX}. It could
be regarded as a section of $\wedge^3(L^*)^*\cong \wedge^3 L$. It is
called the ``Nijenhuis operator" in \cite[Proposition 3.16]{Marco}.
It is denoted by $\nij$. Its relation with Jacobi identity is
explicitly given in \cite{Marco}. Since $L^*$ is maximally isotropic
in $(TM\oplus T^*M)_\CC$ with respect to the symmetric paring, the
obstruction for $ \Gamma(L^*)$ being closed with respect to the
Courant bracket is the restriction of $\nij$ on $\Gamma(L^*)$
\cite[Proposition 3.27]{Marco}. To be precise, for any three
sections $v_0, v_1, v_2$ of $\Gamma(L^*)$,
\begin{equation}\label{phi}
\nij(v_0, v_1, v_2)=\frac13\big( \inner{\lbra{v_0}{v_1}}{v_2}
+\inner{\lbra{v_1}{v_2}}{v_0}+\inner{\lbra{v_2}{v_0}}{v_1} \big).
\end{equation}

To compute $\nij$, recall that  $L_F$
is rank-1 and $L=L_{F}\oplus E^{(1,0)}$. Therefore, $\nij$ has two components due to the
decomposition
\[
\wedge^3L=\left(L_F\oplus \wedge^2 E^{(1,0)} \right)\oplus \wedge^3 E^{(1,0)}.
\]
Now assume that $\Gamma(L)$ is closed under the Courant bracket. By conjugation,
\[
\lbra{ \Gamma(E^{(0,1)})}{ \Gamma(E^{(0,1)})}\subseteq
\Gamma(L_{F}\oplus E^{(0,1)})= \Gamma(\oL).
\]
 Since $\oL$ is isotropic,
$\inner{E^{(0,1)}}{L_{F}\oplus E^{(0,1)}}=0$. Therefore, if $v_0,
v_1, v_2$ are all sections of $E^{(0,1)}$, $\nij(v_0, v_1, v_2)=0$.
Hence, up to
permutation $\nij$ is uniquely determined by
\begin{equation}\label{special phi}
\nij(v_0, v_1, \eta)=\frac13\big( \inner{\lbra{v_0}{v_1}}{\eta}
+\inner{\lbra{v_1}{\eta}}{v_0}+\inner{\lbra{\eta}{v_0}}{v_1} \big),
\end{equation}
where $v_0$ and $v_1$ are sections of $E^{(0,1)}$.

\begin{prop}\label{type 2 0} The Nijenhuis operator $\nij$ for a generalized contact
structure ${\cal J}=({F}, \eta, \pi, \theta, \varphi)$ is equal to
\begin{equation}\label{obstruction}
\nij=-\frac12 {F}\wedge (\rho^*d\eta)^{(2,0)},
\end{equation}
where $(\rho^*d\eta)^{(2,0)}$ is the $\wedge^2E^{(1,0)}$-component
of the pull-back of $d\eta$ via the anchor map $\rho:L^*\to TM$.
\end{prop}
\bproof Suppose that $X$ and $Y$ are sections of $\ker\eta$ and
$\alpha$ and $\beta$ are sections of $\ker{F}$. Let
$v_0=X+\alpha+i\Phi(X+\alpha)$ and $v_1=Y+\beta+i\Phi(Y+\beta)$. In
terms of the components of $\Phi$,
\[
v_0=X+\alpha+i\pi^\sharp\alpha+i\varphi X+i\theta^\flat
X-i\varphi^*\alpha, \quad \mbox{and} \quad \rho(v_0)=X+i\varphi
X+i\pi^\sharp\alpha.
\]
There is a similar expression for $Y+\beta+i\Phi(Y+\beta)$. Note
that the Courant bracket between any 1-forms is equal to zero, and
the space of 1-forms is isotropic with respect to the symmetric
bilinear pairing (\ref{symmetric pairing}). It follows that
\begin{eqnarray*}
&&\nij(v_0, v_1, \eta)\\
&=&\frac13\Big(
\inner{\lbra{\rho(v_0)}{\rho(v_1)}}{\eta}
+\inner{\lbra{\rho(v_1)}{\eta}}{\rho(v_0)}  +\inner{\lbra{\eta}{\rho(v_0)}}{\rho(v_1)}
\Big)\\
&=& \frac16\Big(
\eta(\lbra{\rho(v_0)}{\rho(v_1)}) +\big(\call_{\rho(v_1)}\eta\big)\rho(v_0)
-\big(\call_{\rho(v_0)}\eta\big)\rho(v_1)\Big).
\end{eqnarray*}
Since $\eta(X+i\varphi X+i\pi^\sharp\alpha)=\eta(Y+i\varphi
Y+i\pi^\sharp\beta)=0$, the above is equal to
\begin{eqnarray*}
&&\frac16\Big(
\eta(\lbra{\rho(v_0)}{\rho(v_1)}) +\big(\iota_{\rho(v_1)}d\eta\big)\rho(v_0)
-\big(\iota_{\rho(v_0)}d\eta\big)\rho(v_1)\Big)\\
&=&\frac16\Big(
-d\eta({\rho(v_0)},{\rho(v_1)}) +\big(\iota_{\rho(v_1)}d\eta\big)\rho(v_0)
-\big(\iota_{\rho(v_0)}d\eta\big)\rho(v_1)\Big)\\
&=&-\frac12 d\eta(\rho(v_0), \rho(v_1)).
\end{eqnarray*}
Therefore, $\nij$ is given as claimed.
 \eproof

Note that if $f$ is a function without zero such that
$F'=\frac{1}{f}F$, and $\eta'={f}\eta$, Then on $(\ker \eta\oplus
\ker F)_\CC$, $d\eta' ={f}d\eta$.  Therefore, the equality in
(\ref{obstruction}) is independent of choice of representative
tensors within a given generalized contact structure.

\

Suppose that $L$ and $L^*$ are both Lie algebroids. Let $d_L$ be the
Lie algebroid differential associated to the bracket on $L$. It acts
on the space of sections of $\wedge^k(L^*)$. Similarly, we have a
differential $d_{L^*}$ associated to the Lie algebroid structure of
$L^*$. It acts on sections of $\wedge^kL$. Since both $L$ and $L^*$
inherit the bracket from the Courant bracket on $(TM\oplus
T^*M)_\CC$, and they are dual to each other with respect to the
symmetric pairing (\ref{symmetric pairing}), they together naturally
form a Lie bialgebroid \cite{LWX}.
%i.e. for any sections $\om_1,
%\om_2$ of $L^*$ and any sections $\ell_1, \ell_2$ of $L$,
%\begin{equation}
% d_L\lbra{\om_1}{\om_2}=\lbra{d_L\om_1}{\om_2}+\lbra{\om_1}{d_L\om_2}.
% \quad
%d_{L^*}\lbra{\ell_1}{\ell_2}=\lbra{d_{L^*}\ell_1}{\ell_2}+\lbra{\ell_1}{d_{L^*}\ell_2}.
%\end{equation}
To summarize our discussion so far, we have the following theorem.

\begin{thm}\label{bi} Let ${\cal J}=({F}, \eta, \pi, \theta,
\varphi)$ represent an (integrable) generalized contact structure.
The pair $L$ and $L^*$ forms a Lie bialgebroid if and only if
$d\eta$ is type (1,1) with respect to the map $\Phi$ on
$(\ker\eta\oplus \ker F)_\CC$.
\end{thm}
\bproof Since $\eta$ is a real 1-form, $d\eta$ is a real 2-form.
Therefore, $(\rho^*d\eta)^{(2,0)}$ is the complex conjugation of
$(\rho^*d\eta)^{(0,2)}$. Therefore, $(\rho^*d\eta)^{(0,2)}=0$ if and
only if $(\rho^*d\eta)^{(2,0)}=0$. \eproof

 The above analysis indicates a special class of objects among generalized
 contact structures.

 \begin{dfn}\label{odd} An almost generalized contact structure
  is called a strong generalized contact structure
   if both $\Gamma(L)$ and $\Gamma(L^*)$ are closed under the
Courant bracket.
\end{dfn}

\subsection{Integrability of  the associated complex subbundles}
Suppose that ${\cal J}=({F}, \eta, \pi, \theta, \varphi)$ represents
a strong generalized contact structure. By complex conjugation, the
closedness of $\Gamma(L^*)$ is equivalent to the closedness of
$\Gamma(\oL^*)$. Since $L=L_F\oplus E^{(1,0)}$ and
$\oL^*=L_\eta\oplus E^{(1,0)}$,
\begin{equation}
\lbra{\Gamma(E^{(1,0)})}{\Gamma(E^{(1,0)})}\subseteq
\Gamma(L_F\oplus E^{(1,0)})\bigcap \Gamma(L_\eta\oplus E^{(1,0)}).
\end{equation}
This inclusion implies
\begin{equation}\label{d eta}
\lbra{\Gamma(E^{(1,0)})}{\Gamma(E^{(1,0)})}\subseteq
\Gamma(E^{(1,0)})
\end{equation}
and the corresponding statement with a complex conjugation.

 With respect to the
symmetric non-degenerate bilinear pairing (\ref{symmetric pairing}),
the dual of $E^{(1,0)}$ is its conjugate bundle $E^{(0,1)}$. In this
section, we focus on the structures of these two bundles.

Our issue now is  whether the pair $E^{(1,0)}$ and $E^{(0,1)}$ forms
a Lie bialgebroid. Since both bundles are Lie algebroids, the only
point for concern is whether there is a natural compatibility
between Lie algebroid differentials and the Courant brackets.

By natural compatibility, we mean to treat the bundles $E^{(1,0)}$
and $E^{(0,1)}$ as subbundles of $(TM\oplus T^*M)_\CC$ with the
Courant bracket (\ref{courant bracket}). While $(TM\oplus T^*M)_\CC$
is a Courant algebroid, the direct sum $E^{(1,0)}\oplus E^{(0,1)}$
may fail to be one because the bracket between sections of
$E^{(1,0)}$ and $E^{(0,1)}$ with respect to the Courant bracket on
$(TM\oplus T^*M)_\mathbb{C}$ may not be a section of
$E^{(1,0)}\oplus E^{(0,1)}$.

 Let $\omega$ be a section of $E^{(1,0)}$ and
$\overline\sigma$ be a section of $E^{(0,1)}$. Suppose that the pair
$E^{(1,0)}$ and $E^{(0,1)}$ forms a Lie bialgebroid. Then
$\lbra{\omega}{\overline\sigma}$ is a section of $E^{(1,0)}\oplus
E^{(0,1)}$. In particular,
$\eta(\rho\lbra{\omega}{{\overline\sigma}})=\eta([\rho(\omega)
,\rho({\overline\sigma})])=0$. By definitions of $E^{(1,0)}$ and
$E^{(0,1)}$, $\eta(\rho(\omega))=0$ and
$\eta(\rho({\overline\sigma}))=0$. Therefore, \[ d\eta(\rho(\omega),
\rho({\overline\sigma}))=-\eta([\rho(\omega)
,\rho({\overline\sigma})])=0.
\]
 In other words, $(\rho^*d\eta)^{(1,1)}$
vanishes identically on $\ker\eta$. As we assume that $d\eta$ is
type $(1,1)$ in the first place, it follows that $d\eta$ vanishes
identically on $\ker\eta$.

Conversely, let $d_{E}$ be the Lie algebroid differential for
$E^{(1,0)}$, and $d_{{\overline E}}$ for $E^{(0,1)}$. The
differential $d_E$ is the composition of an inclusion map, the
differential $d_L$ and a projection. More precisely, given
$\Gamma\left(\wedge ^kL^*\right)= \Gamma \left(\wedge^k(L_\eta\oplus
E^{(0,1)})\right)$ for all $k$, the differential $d_E$ is given by
\[
d_E: \Gamma \left( \wedge^kE^{(0,1)} \right)\hookrightarrow  \Gamma \left(\wedge
^kL^* \right)\stackrel{d_L}\longrightarrow \Gamma \left(\wedge^{k+1}L^*\right)
\stackrel{p}\rightarrow \Gamma \left(\wedge^{k+1}E^{(0,1)}\right),
\]
where the map
\[
p: \Gamma \left(\wedge^{k+1}L^*\right)= \Gamma \left(
\left(L_\eta\otimes \wedge^{k}E^{(0,1)}\right)\oplus
\left(\wedge^{k+1}E^{(0,1)}\right) \right) \rightarrow
\Gamma(\wedge^{k+1}E^{(0,1)})
\]
is a natural projection. i.e.
$d_E{\overline\alpha}=p(d_L\overline\alpha)$ for each
$\overline\alpha$ section of $\wedge^kE^{(0,1)}$.

To check whether the pair $E^{(1,0)}$ and $E^{(0,1)}$ forms a Lie
bialgebroid, we need to verify if
\begin{equation}\label{bi E}
d_{\oE}\lbra{\om_1}{\om_2}=\lbra{d_{\oE}\om_1}{\om_2}+\lbra{\om_1}{d_{\oE}\om_2}
\end{equation}
for any pair of sections $\om_1$ and $\om_2$ of the bundle
$E^{(1,0)}$. Since the pair $L$ and $L^*$ forms a Lie bialgebroid,
\begin{equation}
d_{L^*}\lbra{\om_1}{\om_2}=\lbra{d_{L^*}\om_1}{\om_2}+\lbra{\om_1}{d_{L^*}\om_2}.
\end{equation}
When $\Omega$ is a section of
$
\wedge^2 L=\left(L_F\otimes E^{(1,0)}\right) \oplus \wedge^2E^{(1,0)},
$
let $\Omega^F$ be the first component of $\Omega$ in this decomposition.
Then the above identity becomes
\begin{eqnarray}
&&(d_{L^*}\lbra{\om_1}{\om_2})^F+d_{\oE}\lbra{\om_1}{\om_2} \nonumber\\
&=&\lbra{(d_{L^*}\om_1)^F}{\om_2}+\lbra{d_{\oE}\om_1}{\om_2} +
\lbra{\om_1}{(d_{L^*}\om_2)^F}+\lbra{\om_1}{d_{\oE}\om_2}. \label{E
and L}
\end{eqnarray}
To calculate $(d_{L^*}\om)^F$ for any section $\om$ of $E^{(1,0)}$,
let $\overline\sigma$ be any section of $E^{(0,1)}$. Then by
definition,
\[
 (d_{L^*}\om)(\eta, {\overline\sigma})=\rho(\eta)\inner{\om}{\overline\sigma}-\rho({\overline\sigma})\inner{\om}{\eta}-
\inner{\om}{\lbra{\eta}{{\overline\sigma}}}.
\]
Since $\rho(\eta)=0$ and $\om$ is a section of $\ker\eta\oplus \ker
F$, the above is reduced to
\begin{eqnarray*}
&&\inner{\om}{\cl_{\rho({\overline\sigma})}\eta}=\inner{\rho(\om)}{\cl_{\rho({\overline\sigma})}\eta}
=\inner{\rho(\om)}{\iota_{\rho({\overline\sigma})}d\eta}
=\frac12d\eta(\rho({\overline\sigma}), \rho(\om))\\
&=&-\frac12(\rho^*d\eta)(\omega,
{\overline\sigma})=-\frac12(\iota_\om\rho^*d\eta)({\overline\sigma}).
\end{eqnarray*}
It follows that as a section of $L_F\otimes E^{(1,0)}\subset
\wedge^2 L$,
\begin{equation}
(d_{L^*}\om)^F=-\frac12 F\wedge(\iota_\om\rho^*d\eta)^{(1,0)}.
\end{equation}
Suppose that $(\rho^*d\eta)^{(1,1)}=0$, then
$d\eta(\rho({\overline\sigma}), \rho(\om))=0$. It means that
$(d_{L^*}\om)(\eta, {\overline\sigma})=0$ for all section $\omega$
of $E^{1,0}$ and $\overline\sigma$ of $E^{0,1}$. Therefore, Identity
(\ref{E and L}) is equivalent to Identity (\ref{bi E}). It means
that the pair $E^{(1,0)}$ and $E^{(0,1)}$ forms a Lie bialgebroid.

Note that we assume that $(\rho^*d\eta)^{(2,0)}=0$ in the first
place, the assumption $(\rho^*d\eta)^{(1,1)}=0$ is equivalent to
$\rho^*d\eta=0$ on $\ker \eta$. The next proposition follows.

\begin{prop}\label{E bialgebroid}
Let ${\cal J}=({F}, \eta, \pi, \theta, \varphi)$ represent a strong
generalized contact structure. The pair $E^{(1,0)}$ and $E^{(0,1)}$
with the induced Courant bracket is a Lie bialgebroid if and only if
$d\eta$ vanishes identically on $\ker\eta$.
\end{prop}

\section{Some Classical Geometry in Odd Dimensions}\label{classical}

\subsection{Contact structures}\label{classical contact}
 Suppose that $M$ is a $(2n+1)$-dimensional
manifold with a 1-form $\eta$ such that $\eta\wedge (d\eta)^n$ is
non-zero everywhere, then the 1-form $\eta$ is a contact 1-form.

To make an almost generalized contact structure associated to the
contact 1-form $\eta$, let $\theta=d\eta$. Then the map
\begin{equation}\label{flat}
\flat(X):=\iota_X\theta-\eta(X)\eta
\end{equation}
is an isomorphism from the tangent bundle to the cotangent bundle. In particular, there is a
unique vector field ${F}$ such that $\iota_{F}\eta=1$  and
$\iota_{F} d\eta=0.$  This vector field is known as the Reeb field
of the contact form $\eta$. Define a bivector field $\pi$ by
\begin{equation}\label{pi for cosymplectic}
\pi(\alpha, \beta):=\theta(\flat^{-1}(\alpha), \flat^{-1}(\beta)).
\end{equation}
Choose $\varphi=0$, and
\begin{equation}\label{Phi for cosymplectic}
\Phi=\left(
\begin{array}{cc}
0 & \pi^\sharp \\
\theta^\flat & 0
\end{array}
\right).
\end{equation}
Then, the map $\Phi$, the Reeb field $F$ and the contact form $\eta$
define an almost generalized contact structure.

As the differential forms $\eta$ and $\theta$ are invariant with
respect to the Reeb field $F$, the map $\Phi$ is also invariant.
Therefore,
\[
\lf\eta=0, \quad \lf \theta=0, \quad \lf \Phi=0.
\]
Next we examine the properties of the associated bundles $L$ and
$L^*$.

By Darboux Theorem \cite{GS}, in a neighborhood of any point on $M$,
there exist local coordinates $(x_1, y_1, \dots, x_n, y_n, z)$ such
that
\begin{equation}
 \eta=dz-\sum_{j=1}^n y_jdx_j.
 \end{equation}
The Reeb field is naturally $F=\frac{\partial}{\partial z}$, and $
\theta=d\eta=\sum_{j=1}^n dx_j\wedge dy_j. $ Let
\[
X_j=\frac{\partial}{\partial x_j}+y_j\frac{\partial}{\partial z},
\quad Y_j=\frac{\partial}{\partial y_j}.
\]
Then $\{X_j, Y_j, F\}$ forms a moving frame on the given coordinate,
and $\{dx_j, dy_j, \eta\}$ forms a co-frame.
By construction (\ref{flat}),
\[
\flat(X_j)=dy_j, \quad \flat(Y_j)=-dx_j, \quad \flat(F)=-\eta.
\]
By (\ref{pi for cosymplectic}),
\begin{equation}
\pi=\sum_{j=1}^n X_j\wedge Y_j.
\end{equation}
It follows that $\Phi(\eta)=0$,
\begin{eqnarray*}
\Phi(X_j)=\theta^\flat(X_j)=dy_j, &\quad \mbox{and} \quad
 & \Phi(Y_j)=\theta^\flat(Y_j)=-dx_j. \\
\Phi(dx_j)=\pi^\sharp(dx_j)=Y_j, &\quad \mbox{and} \quad &
\Phi(dy_j)=\pi^\sharp(dy_j)=-X_j.
\end{eqnarray*}
Then a local frame for $E^{(1,0)}$ is $\{X_j-idy_j,  Y_j+idx_j\}.$ A
local frame for $E^{(0,1)}$ is $\{ X_j+idy_j, Y_j-idx_j\}.$ Since
\begin{eqnarray}
&\lbra{F}{X_j-idy_j}=0, \quad \lbra{F}{Y_j+idx_j}=0,& \\
 &\lbra{X_j-idy_j}{Y_j+idx_j}=\lbra{X_j}{Y_j}=-F,&
\end{eqnarray}
the spaces of sections of the bundles $L=L_F\oplus E^{(1,0)}$ and $\oL=L_F\oplus E^{(0,1)}$
are closed under the Courant bracket. It explicitly shows
that $L$ and $\oL$ are Lie algebroids. On the other hand,
\[
\lbra{X_j-idy_j}{\eta}=\iota_{X_j}d\eta=dy_j, \quad
\lbra{Y_j+idx_j}{\eta}=\iota_{Y_j}d\eta=-dx_j.
\]
Therefore, $\Gamma(L^*)=\Gamma(L_\eta\oplus E^{(1,0)})$ and $\Gamma(\oL^*)=\Gamma(L_\eta\oplus
E^{(0,1)})$ are not closed under the Courant bracket.

 As the obstruction to closedness is $F\wedge
(\rho^*d\eta)^{2,0}$, we could also find it through the
type-decomposition of $\rho^*d\eta$. Given $\theta =
d\eta=\sum_{j=1}^n dx_j\wedge dy_j$ and the map $\Phi$ above, it is
straightforward to find that
%\begin{eqnarray*}
%\theta &=& d\eta=\sum_{j=1}^n dx_j\wedge dy_j\\
%&=& \sum_{j=1}^n\Big( \big( \frac12(dx_j-iY_j)+\frac12(dx_j+iY_j)
%\big) \wedge \big( \frac12(dy_j+iX_j)+\frac12(dy_j-iX_j) \big)
% \Big)
%\end{eqnarray*}
%Therefore,
\begin{eqnarray*}
(\rho^*d\eta)^{2,0} &=&\frac14\sum_{j=1}^n(dx_j-iY_j)\wedge (dy_j+iX_j),\\
(\rho^*d\eta)^{0,2} &=&\frac14\sum_{j=1}^n(dx_j+iY_j)\wedge (dy_j-iX_j),\\
(\rho^*d\eta)^{1,1} &=& \frac12 \sum_{j=1}^n(dx_j\wedge
dy_j+X_j\wedge Y_j)=\frac12(d\eta+\pi).
\end{eqnarray*}
In particular, the obstruction for closedness of $\Gamma(L^*)$ does not
vanish anywhere.

\begin{prop}\label{contact} Let ${\cal J}=({F}, \eta, \pi, \theta,
\varphi)$ represent a generalized contact structure associated to a
classical contact 1-form $\eta$. Then the corresponding bundles $L$
and $\oL$ are Dirac structures. The bundles $L^*$ and $\oL^*$ are
never Dirac structures. In particular, the pair $L$ and $L^*$ is not
a Lie bialgebroid.
\end{prop}

\subsection{Almost cosymplectic structures}\label{almost cosym}

An almost cosymplectic structure consists of   a 1-form $\eta$ and a
2-from $\theta$ such that $\eta\wedge\theta^n\neq 0$ at every point
of the manifold. Given this condition,  the map formally given  in
(\ref{flat}) is again an isomorphism. Therefore, there exists a
unique vector field $F$ such that $\eta(F)=1$ and $\theta(F)=0$.
 These
tensors determine an almost generalized contact structure by the
matrix $\Phi$. It is  formally given in (\ref{Phi for cosymplectic}).

If both $\eta$ and $\theta$ are closed,  we address the pair $(\eta, \theta)$ a
cosymplectic structure without qualification.
Next, we investigate integrability of the generalized almost contact
structure associated to a cosymplectic structure $(\eta, \theta)$.
Since $\iota_F\theta=0$, for any section $X$ of $\ker \eta$,
$
\lbra{F}{X-i\iota_{X}\theta}=[F, X]-i\cl_F\iota_X\theta.
$
Since $\theta$ is closed and $\iota_F\theta=0$, $\cl_F\theta=0$. As
$\cl_F\iota_X\theta-\iota_X\cl_F\theta=\iota_{[F, X]}\theta$, it
follows that
\[
\lbra{F}{X-i\iota_{X}\theta}=[F, X]-i\iota_{[F,X]}\theta.
\]
If $X$ and $Y$ are sections of $\ker\eta$,
\begin{eqnarray*}
&&\lbra{X-i\iota_X\theta}{Y-i\iota_Y\theta}
=[X,Y]-i(\cl_X\iota_Y\theta-\cl_Y\iota_X\theta)+id(\iota_X\iota_Y\theta)\\
&=&[X,Y]-i\iota_{[X,Y]}\theta-i(\iota_Y\cl_X\theta-\cl_Y\iota_X\theta-d\iota_X\iota_Y\theta).
\end{eqnarray*}
It is equal to $[X,Y]-i\iota_{[X,Y]}\theta$ due to $d\theta=0$ and
the identity $\cl_X=d\circ \iota_X+\iota_X\circ d$.

Through the isomorphism $\flat$,  the computation above also shows
that for any sections $\alpha$ and $\beta$ of $\ker F$,
$\lbra{\alpha-i\iota_\alpha\pi}{\beta-i\iota_\beta\pi}$ is a section
of $E^{(1,0)}$. Similarly,
$\lbra{X-i\iota_X\theta}{\beta-i\iota_\beta\pi}$ is a section of
$E^{(1,0)}$ whenever both ${X-i\iota_X\theta}$ and
${\beta-i\iota_\beta\pi}$ are. Therefore,  $\Gamma(E^{(1,0)})$ is
closed under the Courant bracket.

It follows that $\Gamma(L)$ is closed under the Courant bracket. In
addition, since $d\eta=0$, by Theorem \ref{bi}, Definition
\ref{odd}, and Proposition \ref{pi for cosymplectic}, we have the
following observation.

\begin{prop}\label{L and E} If $\cal J$ represents a generalized
almost contact structure associated to a classical cosymplectic
structure, then it is a strong generalized contact structure.
Moreover, the pairs of bundles $(L, L^*)$ and $(E^{(1,0)},
E^{(0,1)})$ with respect to the induced Courant bracket are both Lie
bialgebroids.
\end{prop}

\subsection{Almost contact structures}\label{almost contact}
 Suppose that $M$ is a
(2n+1)-dimensional manifold with a vector field ${F}$, a 1-form
$\eta$ and a type (1,1)-tensor $\varphi$ satisfying
\begin{equation}\label{def almost contact}
\varphi^2=-\mathbb{I}+\eta\otimes {F} \quad \mbox{ and } \quad \eta({F})=1,
\end{equation}
then the triple $(\varphi, {F}, \eta)$ is a called an \it almost
contact structure \rm\cite{Sasaki1}.

Associated to any almost contact structure, we have an almost
generalized contact structure by setting
\begin{equation}\label{Phi for normal contact}
\Phi= \left(
\begin{array}{cc}
\varphi  & 0 \\
0 & - \varphi^*
\end{array}
\right)
\end{equation}
with the given vector field $F$ and 1-form $\eta$.

An almost contact structure  is a ``normal almost contact structure"
\cite{Blair} if
\begin{equation}\label{normal contact}
 {\cal N}_{\varphi}=-{F}\otimes d\eta, \quad \lf \varphi=0 \quad \mbox{ and } \quad \lf
 \eta=0,
 \end{equation}
where by definition,
\begin{equation}
{\cal N}_{\varphi}(X, Y)= [ \varphi X, \varphi Y ] + \varphi^2 [ X,
Y ]
 - \varphi([ \varphi X, Y ] + [ X, \varphi Y ])
 \end{equation}
 for any vector fields $X$ and $Y$.
Note that equations (\ref{normal contact}) imply that  if $s$ is a
section of $E^{(1,0)}$, then $\lbra{F}{s}$ is again a section of
$E^{(1,0)}$.

Since ${\cal N}_\varphi=- F\otimes d\eta$, for any vector fields $X$
and $Y$
\begin{equation}\label{nij phi}
-d\eta(X,Y)F=[\varphi X, \varphi Y]+\varphi^2[X,Y]-\varphi([\varphi
X, Y]+[X, \varphi Y]).
\end{equation}
In particular, this identity holds when the vector fields are
sections of the bundle $\ker\eta$. In such case, applying $\eta$ on
both sides of this identity, we find that
\begin{equation}
\eta([\varphi X,\varphi Y])=\eta({\cal N}_\varphi(X,
Y))=-d\eta(X,Y).
\end{equation}
As $\varphi X$ and $\varphi Y$ are also sections of $\ker\eta$, the
above identity implies that
\begin{equation}\label{type 11 formula}
 d\eta(\varphi X, \varphi Y)=d\eta (X, Y)
\end{equation} for any sections $X$ and $Y$ in $\ker\eta$.
Therefore, the restriction of $d\eta$ on
$\ker \eta$ is type-(1,1) with respect to $\Phi$.  For future reference, we highlight
this observation.
\begin{lemma}\label{type 11 lemma} Suppose that $(F, \eta, \varphi)$ is a normal almost
contact
structure, then $\rho^*d\eta$ is a section of $E^{(1,0)}\otimes E^{(0,1)}$.
\end{lemma}

Now, for any sections $X$ and $Y$ of $\ker\eta$, due to the first
identity in (\ref{def almost contact})
\begin{eqnarray}
&&\lbra{X-i\varphi X}{Y-i\varphi Y}\label{intermediate}\\
&=&[X,Y]-[\varphi X,\varphi Y]-i([\varphi X, Y]+[X, \varphi Y])\nonumber\\
&=&[X,Y]+\varphi^2[\varphi X,\varphi Y]-i([\varphi X, Y]+[X, \varphi
Y])-\eta([\varphi X,\varphi Y])F.\nonumber
\end{eqnarray}
Since $\varphi X$ and $\varphi Y$ are sections of $\ker\eta$,
$\eta([\varphi X,\varphi Y])=-d\eta(\varphi X,\varphi Y)$. Applying
formula (\ref{nij phi}) to the pair of vector fields $\varphi X$ and
$\varphi Y$, and observing that $\varphi^2X=-X$, $\varphi^2Y=-Y$, we
get
\[
-d\eta(\varphi X,\varphi Y)F=[ X, Y]+\varphi^2[\varphi X, \varphi
Y]+\varphi([ X, \varphi Y]+[\varphi X, Y]).
\]
Given (\ref{def almost contact}) and (\ref{type 11 formula}),
(\ref{intermediate}) is equal to
\begin{eqnarray*}
&&-\varphi ([\varphi X, Y]+[X, \varphi Y])-i([\varphi X, Y]+[X,
\varphi Y])\\
&=&-\varphi ([\varphi X, Y]+[X, \varphi Y])+i\varphi^2([\varphi X,
Y]+[X, \varphi Y])\\
&&-i\eta([\varphi X, Y]+[X, \varphi Y])F\\
&=&-\varphi ([\varphi X, Y]+[X, \varphi Y])+i\varphi^2([\varphi X,
Y]+[X, \varphi Y]).
\end{eqnarray*}
Therefore,
\begin{eqnarray*}
&&\lbra{X-i\varphi X}{Y-i\varphi Y}\\
&=&-\varphi ([\varphi X, Y]+[X, \varphi Y])+i\varphi^2([\varphi X,
Y]+[X, \varphi Y]).
\end{eqnarray*}
Since $-\varphi ([\varphi X, Y]+[X, \varphi Y])$ is a section of
$\ker \eta$, the above tensor is a section of $E^{(1,0)}$.

Next, suppose that $X$ is a section of $\ker \eta$ and $\beta$ is a
section of $\ker F$.  By definition
\begin{eqnarray*}
&&\lbra{X-i\varphi X}{\beta+i\varphi^* \beta}\\
&=& \cl_{(X-i\varphi X)}{(\beta+i\varphi^* \beta)}-\frac12
d\iota_{(X-i\varphi X)}(\beta+i\varphi^* \beta)\\
&=&\cl_{(X-i\varphi X)}{(\beta+i\varphi^* \beta)}\\
&=&\cl_X\beta+\cl_{(\varphi X)}(\varphi^*\beta)+i(\cl_X(\varphi^*
\beta)-\cl_{(\varphi X)}\beta).
\end{eqnarray*}
Evaluating the real part of the above expression on the Reeb field,
with standard tensor calculus and (\ref{def almost contact}), we
find that it is equal to $\eta([F, X])\beta(F)$. Since $\beta$ is a
section of $\ker F$, the real part of the above expression is a
section of $\ker F$. Next, due to transpose of the first formula in
(\ref{def almost contact}),
\begin{eqnarray*}
&& \varphi^*(\cl_X\beta+\cl_{\varphi X}(\varphi^*\beta))\\
%&=&\varphi^*(\cl_X\beta)+\varphi^*(\cl_{(\varphi
%X)}(\varphi^*\beta))\\
&=&\cl_X(\varphi^*\beta)-(\cl_X\varphi)^*\beta+\cl_{(\varphi
X)}((\varphi^*)^2\beta)-(\cl_{(\varphi
X)}\varphi)^*(\varphi^*\beta)\\
&=&\cl_X(\varphi^*\beta)-(\cl_X\varphi)^*\beta+\cl_{(\varphi
X)}(-\beta+\beta(F)\eta)-(\cl_{(\varphi
X)}\varphi)^*(\varphi^*\beta)\\
&=&\cl_X(\varphi^*\beta)-\cl_{(\varphi
X)}\beta-(\cl_X\varphi)^*\beta-(\cl_{(\varphi
X)}\varphi)^*(\varphi^*\beta).
\end{eqnarray*}
We claim that $(\cl_X\varphi)^*\beta+(\cl_{(\varphi
X)}\varphi)^*(\varphi^*\beta)=0$. To verify, let $A$ be any vector
field,
\begin{eqnarray*}
&&\big( (\cl_X\varphi)^*\beta\big)A+\big((\cl_{(\varphi
X)}\varphi)^*(\varphi^*\beta)\big)A\\
&=&\beta\big( (\cl_X\varphi)A +\varphi\circ (\cl_{(\varphi
X)}\varphi)A \big)\\
%&=&\beta \big( \cl_X(\varphi
%A)-\varphi(\cl_XA)+\varphi(\cl_{(\varphi X)}(\varphi
%A)-\varphi(\cl_{\varphi X}A)) \big)\\
&=&\beta\big( [X, \varphi A]-\varphi [X, A]+\varphi [\varphi X,
\varphi A]-\varphi^2 [\varphi X, A]\big)\\
&=&\beta\big( [X, \varphi A]-\varphi [X, A]+\varphi [\varphi X,
\varphi A]+ [\varphi X, A]-\eta([\varphi X, A])F \big)\\
%&=&\beta\big( [X, \varphi A]-\varphi [X, A]+\varphi [\varphi X,
%\varphi A]+ [\varphi X, A] \big)\\
&=&\beta \big( \varphi {\cal N}_{\varphi}(X, A)\big).
\end{eqnarray*}
By (\ref{normal contact}), ${\cal N}_{\varphi}(X, A)=-d\eta(X, A)
F$. Since $\varphi(F)=0$, $\varphi {\cal N}_{\varphi}(X, A)=0$.

Since the Courant bracket between two 1-forms is always equal to
zero, we could now conclude that the Courant bracket between two
sections of $E^{(1,0)}$ is again a section of $E^{(1,0)}$. Since the
bundle $E^{(1,0)}$ is $F$-invariant, the bundle $L$ is closed with
respect to the Courant bracket.

Finally, formula (\ref{type 11 formula}) shows that
$(\rho^*d\eta)^{(2,0)}=0$. Therefore $L^*$ is closed with respect to the
Courant bracket. By Theorem \ref{bi} we have the following.

\begin{prop}\label{normal and generalized} If $\cal J$ represents a
generalized almost contact structure associated to a classical
normal almost contact structure on an odd-dimensional manifold $M$,
then it is a strong generalized contact structure.
\end{prop}

\section{Examples of Strong Generalized Contact Structures}

\subsection{Structures on SU(2)}\label{SU2}
On the Lie algebra $\su(2)$, choose a basis $X_1, X_2, X_3$ and dual
basis $\sigma^1, \sigma^2, \sigma^3$ such that
\begin{equation}
[X_1, X_2]=-X_3, \quad d\sigma^1=\sigma^2\wedge \sigma^3,
\end{equation}
and cyclic permutations of the indices $\{1,2,3\}$.

\subsubsection{Normal contact structures on SU(2)}
To construct a classical normal almost contact structure, one simply
takes
\begin{equation}
\eta=\sigma^3, \quad F=X_3, \quad \varphi=X_2\otimes
\sigma^1-X_1\otimes \sigma^2.
\end{equation}
Then
\begin{equation}
\varphi^*=\varphi=-\sigma^2\otimes X_1+\sigma^1\otimes X_2.
\end{equation}
Therefore,
\begin{eqnarray*}
&\Phi(X_1)=\varphi(X_1)=X_2, \quad \Phi(X_2)=\varphi(X_2)=-X_1,&\\
&\Phi(\sigma^1)=-\varphi^*(\sigma^1)=\sigma^2, \quad
\Phi(\sigma^2)=-\varphi^*(\sigma^2)=-\sigma^1.&
\end{eqnarray*}
The bundle $L$ and $L^*$ are globally trivialized. As modules over
the space of smooth functions,
\begin{eqnarray*}
&\Gamma(L)=\Gamma(L_F\oplus E^{1,0})=\vspan{X_3, \frac{1}{\sqrt 2}(X_1-iX_2),
\frac{1}{\sqrt 2}(\sigma^1-i\sigma^2)},& \\
&\Gamma(L^*)=\Gamma(L_\eta\oplus E^{0,1})=\vspan{\sigma^3, \frac{1}{\sqrt
2}(X_1+iX_2), \frac{1}{\sqrt 2}(\sigma^1+i\sigma^2)}.&
\end{eqnarray*}
It is now an elementary computation to verify that the structure
equations for Lie algebroids $L$ and $L^*$ are respectively given by
\begin{eqnarray*}
&\lbra{X_3}{\frac{1}{\sqrt 2}(X_1-iX_2)}=-\frac{i}{\sqrt
2}(X_1-iX_2),& \\
&\hspace{1in}\lbra{X_3}{\frac{1}{\sqrt
2}(\sigma^1-i\sigma^2)}=-\frac{i}{\sqrt 2}(\sigma^1-i\sigma^2),&\\
&\hspace{2in}\lbra{\sigma^3}{\frac{1}{\sqrt 2}(X_1+iX_2)}
=\frac{i}{\sqrt 2}(\sigma^1+i\sigma^2).&
\end{eqnarray*}
On the other hand, we have
\begin{equation}
\lbra{X_1-iX_2}{X_1+iX_2}=-2iX_3.
\end{equation}
It demonstrates that $\Gamma(E^{1,0}\oplus E^{0,1})$
 is not closed under the Courant
bracket. In other words, with respect to the induced Courant
bracket, $E^{1,0}\oplus E^{0,1}$ is not a Courant algebroid
\cite{LWX}. It follows that the pair $E^{1,0}$ and $E^{0,1}$, with
respect to the induced Courant bracket, does not form a Lie
bialgebroid. This example demonstrates that Proposition \ref{L and
E} for cosymplectic structures could not be extended to normal
almost contact structures, or  strong generalized contact structures
in general.

\subsubsection{Contact structures on SU(2)}
An obvious contact
structure on $\SU(2)$ is given by $\eta=\sigma^3$. In such case,
\begin{equation}
F=X_3, \quad \theta=d\sigma^3=\sigma^1\wedge\sigma^2, \quad
\pi=X_1\wedge X_2.
\end{equation}
With $\varphi=0$, the restriction of $\Phi$ on $\ker\sigma^3\oplus
\ker X_3$ is determined by
\begin{equation}
\Phi(X_1)=\sigma^2, \quad \Phi(X_2)=-\sigma^1, \quad
\Phi(\sigma^1)=X_2, \quad \Phi(\sigma^2)=-X_1.
\end{equation}
Therefore,
\begin{equation}
L=\vspan{X_3, X_1-i\sigma^2, X_2+i\sigma^1}, \quad
L^*=\vspan{\sigma^3, X_1+i\sigma^2, X_2-i\sigma^1}.
\end{equation}
Taking the Courant brackets, we find that
\begin{eqnarray*}
&\lbra{X_3}{X_1-i\sigma^2}=-(X_2+i\sigma^1), \quad
\lbra{X_3}{X_2+i\sigma^1}=X_1-i\sigma^2,&\\
&\lbra{X_1-i\sigma^2}{X_2+i\sigma^1}=-X_3=\lbra{X_1+i\sigma^2}{X_2-i\sigma^1},& \\
&\lbra{\sigma^3}{X_1+i\sigma^2}=-\sigma^2, \quad
\lbra{\sigma^3}{X_2-i\sigma^1}=\sigma^1.
\end{eqnarray*}
This example reaffirms that $L$ forms a Lie algebroid while $L^*$
fails to be one.

\subsection{Structures on the 3-dimensional Heisenberg group}\label{H3}
On the three-dimensional Heisenberg group $H_3$, we choose a basis
$\{X_1, X_2, X_3\}$ for its algebra $\lie h_3$ so that $[X_1,
X_2]=-X_3$. Let $\{\alpha^1, \alpha^2, \alpha^3\}$ be a dual frame.
Then $d\alpha^3=\alpha^1\wedge \alpha^2$.

\subsubsection{Cosymplectic structure on $H_3$}
For any real numbers $a$ and $b$,
choose
\begin{equation}\label{ab}
\eta=\alpha^1 \quad \mbox{and} \quad
\theta=\alpha^2\wedge\alpha^3+a\alpha^1\wedge\alpha^2+b\alpha^1\wedge\alpha^3.
\end{equation}
They together define a cosymplectic structure. The Reeb field is
$F=X_1-bX_2+aX_3.$ Since \[ \flat(X_1)=a\alpha^2+b\alpha^3-\alpha^1,
\quad \flat(X_2)=\alpha^3-a\alpha^2, \quad
\flat(X_3)=-b\alpha^1-\alpha^1,
\]
$\pi=X_2\wedge X_3$ and $\varphi=0$. Apparently, $\ker F=\langle
\alpha^2+b\alpha^1, \alpha_3-a\alpha^1\rangle$ and $\ker\eta=\langle
X_2, X_3\rangle$. Since
\[
\Phi(X_2)=\alpha^3-a\alpha^1, \quad \Phi(X_3)=-\alpha^2-b\alpha^1,
%,
%\quad \Phi(\alpha^2+b\alpha^1)=X_3, \quad
%\Phi(\alpha^3-a\alpha^1)=-X_2.
\]
we obtain global sections to trivialize the bundles $L$ and $L^*$.
\begin{eqnarray*}
&L=\langle X_1-bX_2+aX_3,
X_2-i\alpha^3+ia\alpha^1, X_3+i\alpha^2+ib\alpha^1\rangle,&\\
&L^*=\langle \alpha^1, X_2+i\alpha^3-ia\alpha^1,
X_3-i\alpha^2-ib\alpha^1\rangle.
\end{eqnarray*}

Since the Courant brackets between $X_3, \alpha^1, \alpha^2$ and any
element among $X_1, X_2, X_3, \alpha^1, \alpha^2, \alpha^3$ are
equal to zero, the restriction of the Courant bracket on $L^*$ is
identically equal to zero. The restriction on $L$ is determined by a
single non-trivial equation, namely
\[
[X_1-bX_2+aX_3,
X_2-i\alpha^3+ia\alpha^1]=-(X_3+i\alpha^2+ib\alpha^1).
\]

\subsubsection{New examples on $H_3$} For $t=rc+irs$ where
$c=\cos\vartheta$ and $s=\sin\vartheta$ for some real number
$\vartheta$,  define
\begin{eqnarray*}
&\varphi_t:=\frac{2rc}{1-r^2}(X_2\otimes \alpha^2+X_3\otimes
\alpha^3),&\\
&\theta_t:= \frac{r^2-2rs+1}{1-r^2}\alpha^2\wedge\alpha^3;\quad
\pi_t= \frac{r^2+2rs+1}{1-r^2}X_2\wedge X_3.&
\end{eqnarray*}
Now as given in (\ref{Phi in general}), define
\[
\Phi_t:= \left(
\begin{array}{cc}
\varphi_t  & \pi_t^\sharp \\
\theta_t^\flat & - \varphi_t^*
\end{array}
\right),
\]
then ${\cal J}_t:=(F, \eta, \pi_t, \theta_t, \varphi_t)$ is a family
of generalized almost contact structures. The corresponding bundles
$L_t$ and its conjugate $\oL_t$ are trivialized.
\begin{eqnarray*}
L_t &=&\langle X_1, (X_2-i\alpha^3)-i\Phi_t(X_2-i\alpha^3),
(X_3+i\alpha^2)-i\Phi_t(X_3+i\alpha^2) \rangle,\\
&=&\langle X_1, (1+rs)X_2+rc\alpha^3-i(1-rs)\alpha^3-ircX_2,\\
&& \qquad \qquad (1+rs)X_3-rc\alpha^2+i(1-rs)\alpha^2-ircX_3 \rangle\\
 L_t^* &=&\langle \alpha^1,
(\alpha^2+iX_3)+i\Phi_t(\alpha^2+iX_3),
(\alpha^3-iX_2)+i\Phi_t(\alpha^3-iX_2)\rangle\\
&=&\langle \alpha^1,
(1-rs)\alpha^2-rcX_3+i(1+rs)X_3-irc\alpha^2,\\
& & \qquad \qquad (1-rs)\alpha^3+rcX_2-i(1+rs)X_2-irc\alpha^3
\rangle.
\end{eqnarray*}
Since the Courant brackets between $X_3, \alpha^1, \alpha^2$ and any
element among $X_1, X_2, X_3, \alpha^1, \alpha^2, \alpha^3$ are
equal to zero, it is straightforward to check that the restriction
of the Courant bracket to $\Gamma(L^*_t)$ is trivial. On
$\Gamma(L_t)$, the sole non-trivial bracket is due to
\begin{eqnarray*}
&&\lbra{X_1}{(1+rs)X_2+rc\alpha^3-i(1-rs)\alpha^3-ircX_2}\\
&=&-\big( (1+rs)X_3-rc\alpha^2+i(1-rs)\alpha^2-ircX_3
 \big).
 \end{eqnarray*}

Therefore, ${\cal J}_t$ is an analytic family of strong generalized
contact structures.

In this family, there are two apparent sub-families, determined by
$|t|^2=r^2<1$ and $|t|^2=r^2>1$.

When $t=0$, we recover the strong generalized contact structure
determined by a cosymplectic structure as given in (\ref{ab}) with
$a=b=0$.

When $r\neq 0$ and $\cos\vartheta\neq 0$, the strong generalized
contact structure is no longer given by a classical cosymplectic
structure. Since the polynomials $r^2-2r\sin\vartheta+1$ and
$r^2+2r\sin\vartheta+1$ do not have zeroes for any $\vartheta$, the
family does not contain any classical almost contact structures
neither.

When $r\to \infty$, we recover the cosymplectic structure with
1-form $\eta=\alpha^1$ and 2-form $\theta_\infty=-\alpha^2\wedge
\alpha^3$.

\subsubsection{Deformation of cosymplectic structures}\label{def of gcx}

Recall Proposition \ref{L and E} that the pair of bundles $E^{1,0}$
and $E^{0,1}$ forms a Lie bialgebroid with respect to the
restriction of the Courant bracket when they are determined by a
classical cosymplectic structure. Let the Lie algebroid differential
for the former to be denoted by $d_E$ and the latter to be denoted
by $d_{\oE}$. Suppose that $\Gamma$ is a section of
$\wedge^2E^{0,1}$. It is also treated as a section of Hom$(E^{1,0},
E^{0,1})$. If it satisfies the Maurer-Cartan equation,
\begin{equation}\label{Maurer-Cartan}
d_E\Gamma+\frac12\lbra{\Gamma}{\Gamma}=0,
\end{equation}
the graph of $\Gamma$ is a deformation of the Lie algebroid
$E^{1,0}$ \cite{LWX}. Denote it by $E^{1,0}_\Gamma$.

Since $E^{0,1}$ is a complex conjugation of $E^{1,0}$, the graph of
$\overline{\Gamma}$ determines a deformation of $E^{0,1}$,
$E^{0,1}_{\overline \Gamma}$. As it is obvious that
$E^{0,1}_{\overline \Gamma}$ is isomorphic to the complex
conjugation of $E^{1,0}_\Gamma$, we obtain a deformation of a
cosymplectic structure through strong generalized contact
structures, with the Reeb field $F$ and the 1-form $\eta$
unperturbed.

In the example of the last section, the restriction of the Courant
bracket on both $E^{1,0}$ and on $E^{0,1}$ are trivial. It follows
that the section
\begin{equation}
\Gamma=(\alpha^2+iX_3)\wedge (\alpha^3-iX_2)
\end{equation}
solves the Maurer-Cartan equation (\ref{Maurer-Cartan}). Therefore,
we obtain deformations. To recover the family of strong generalized
contact structures on $H_3$ in the previous example, one simply
takes $r(\cos\vartheta+i\sin\vartheta)\Gamma$ to generate new
examples.

\section{Examples of Generalized Contact Structures}\label{new}

It is well known that if $\eta$ is a regular contact 1-form on a
compact manifold $M$, then $M$ is a principal circle bundle over a
smooth manifold $N$ such that $\eta$ is a connection 1-form. Here
$N$ is the space of leaves of the foliation of the Reeb field $F$
for the contact 1-form $\eta$. Moreover, there exists a symplectic
form $\omega$ on $N$ such that the curvature form of $\eta$ is given
by $d\eta=-p^*\omega$, where $p:M\to N$ is the quotient map
\cite{BW}. A converse construction of contact structures on any
principal $\SO(2)$-bundle whose characteristic class is a symplectic
form is easily developed through the identity $d\eta=-p^*\omega$. In
this section, we illustrate how these constructions could be done
for generalized contact structures, at least in the case when the
manifolds involved are Lie groups and the geometry are invariant. At
the end, we produce a non-trivial family of generalized contact
structures, with a classical contact 1-form in the family. Thereby,
we demonstrate that classical contact structures have  deformation
in the category of generalized contact structures, and away from
classical objects. It leads to a departure from Gray's theorem that
up to diffeomorphisms, contact structures on compact manifolds do
not have non-trivial deformation among classical contact 1-forms
\cite{Gray}.

\subsection{On central extensions of even-dimensional Lie groups}
\label{contact Lie groups}\label{central extension}

Suppose that $H$ is a real Lie group with an invariant symplectic
form $\omega$. Let $\lie h$ be the Lie algebra of $H$, with Lie
bracket $\bra{-}{-}$. Denote $\lie c$ a one-dimensional real vector
space. Let $F$ be a non-zero vector in $\lie c$. On the space $\lie
g:=\lie h\oplus \lie c$, we next define a new Lie bracket $[-,-]$ on
$\lie g$ as follows. For any $X$ and $Y$ in $\lie h$,
\begin{equation}\label{extended 1}
[X,Y]:=\bra{X}{Y}+\omega(X,Y)F, \qquad \mbox{ and } \qquad [X, F]=0.
\end{equation}
To check that $[-,-]$ is indeed a Lie bracket, one needs only to
check that the Jacobi identity is satisfied by a triple of elements
in $\lie h$. It turns out to be a consequence of  $d\omega=0$ and
$\bra{-}{-}$ satisfying the Jacobi identity. This construction makes
$\lie g$ a central extension of $\lie h$ by $\lie c$.

Elements in $\lie h^*$ are extended to be elements in
$\lie g^*$ by setting their evaluations on $\lie c$ to be equal to zero. Let $\eta$ be the 1-form on $\lie g$ such that
$\eta(X)=0$ for all $X$ in $\lie g$ and $\eta(F)=1$.
%Suppose that $\{X_1, \dots, X_{2n}\}$ is a basis for $\lie h$, and
%$\{\alpha^1, \dots, \alpha^{2n}\}$ is a dual basis. Let the structure equations for $\lie h$ in the dual basis be
%\begin{equation}
%d\alpha^\ell=-\sum_{j<k}^\ell c_{jk}^\ell \alpha^j\wedge \alpha^k.
%\end{equation}
%The structure equations for $\lie g$ are
%\begin{equation}
%d\alpha^\ell=-\sum_{j<k}^\ell c_{jk}^\ell \alpha^j\wedge \alpha^k, \quad \mbox{ with the addition of}  \quad d\eta=-\omega.
%\end{equation}
Next,  for any $X$ in $\lie h$ and any $\alpha$ in $\lie h^*$, we have
\begin{equation}\label{extended 2}
{\cal L}_X\eta=-\iota_X\omega, \quad \mbox{and} \quad  {\cal
L}_F\alpha=0.
\end{equation}

Suppose that
\begin{equation}
\Phi= \left(
\begin{array}{cc}
\varphi  & \pi^\sharp \\
\theta^\flat & - \varphi^*
\end{array}
\right)
\end{equation}
is a generalized complex structure on the Lie group $H$ as given
above. Suppose in addition that all three tensorial components
$\varphi, \theta, \pi$ are left-invariant. Then we treat $\Phi$ as a
real linear map from $\lie h\oplus \lie h^*$, and extend it by zeros
to a linear map from $\lie g\oplus \lie g^*$. It follows that ${\cal
J}=(F, \eta, \pi, \theta, \varphi)$ defines a generalized almost
contact structure on the Lie algebra $\lie g$, and hence as a
left-invariant generalized almost contact structure on the Lie group
$G$, whose algebra is determined by (\ref{extended 1}).

With respect to the notations in Section \ref{bundles} and as far as invariant sections are concerned,
\begin{equation}
\ker \eta=\lie h, \quad \ker F=\lie h^*.
\end{equation}
The spaces of invariant sections of $L$ and $L^*$ are respectively the following finite dimensional complex vector spaces.
\begin{equation}
\lie l=\langle F\rangle_\mathbb{C}\oplus \lie h^{1,0}, \quad \lie l^*=\langle{\eta}\rangle_\mathbb{C}\oplus \lie h^{0,1}.
\end{equation}
Due to the structure equations (\ref{extended 1}) and (\ref{extended
2}),  $ \lbra{\lie l}{\lie l}=\lbra{\lie h^{1,0}}{\lie h^{1,0}}$,
and $\lbra{\lie h^{1,0}}{\lie h^{1,0}}\subseteq \langle
F\rangle_\mathbb{C}\oplus \lie h_\mathbb{C}$. Since $\Phi$ is an
integrable generalized complex structure, the $\lie
h_\mathbb{C}$-component of $\lbra{\lie h^{1,0}}{\lie h^{1,0}}$ is
contained in ${\lie h^{1,0}}$. Therefore,
 $\lbra{\lie l}{\lie l}\subseteq \lie l$.
 From (\ref{extended 1}), we also see that $\lbra{\lie l^*}{\lie
l^*}$ in general is not a subspace of $\lie l^*$. Therefore, we
obtain an invariant generalized contact structure, but not a strong
one.

\subsection{Geometry on four-dimensional Kodaira manifold}
In \cite{Poon}, the first author shows that the complex structure on a primary Kodaira surface could be deformed, within a family
of generalized complex structure, to a symplectic structure. In this section, we briefly recall his construction to establish notations.

A real four-dimensional Kodaira manifold $N$ is a co-compact
quotient of a four-dimensional nilpotent Lie group $H$ \cite{GMPP}.
Let $\{e_1, e_2, e_3, e_4\}$ be a basis of the Lie algebra $\lie h$,
and $\{e^1, \dots, e^4\}$ be the dual basis.
 The sole non-zero structure equation and its dual expression are respectively  given by
\begin{equation}
[e_1, e_2]=e_3 \quad \mbox{and} \quad de^3=-e^1\wedge e^2.
\end{equation}
In particular, the space of invariant closed 2-forms on the Kodaira manifold $N$ is spanned by
\begin{equation}\label{closed 2-forms}
e^1\wedge e^3-e^2\wedge e^4, \quad e^1\wedge e^4+e^2\wedge e^3, \quad e^1\wedge e^3+e^2\wedge e^4, \quad
e^1\wedge e^4-e^2\wedge e^3.
\end{equation}
For any real constants $u_1, v_1, u_2, v_2$ with
$u_1^2+v_1^2-u_2^2-v_2^2\neq 0$,
\begin{eqnarray}
&u_1(e^1\wedge e^3-e^2\wedge e^4)+v_1(e^1\wedge e^4 +e^2\wedge
e^3)&\nonumber\\
&\hspace{1.0in}+u_2(e^1\wedge e^3+e^2\wedge e^4)+v_2(e^1\wedge
e^4-e^2\wedge e^3)& \label{symplectic family}
\end{eqnarray}
is a symplectic form.

On the other hand, the group  $H$ has an invariant integrable complex structure $J$. In terms of the given basis for the Lie algebra $\lie h$,
\[
Je_1=e_2, \quad Je_2=-e_1, \quad Je_3=e_4, \quad Je_4=-e_3.
\]
This complex structure on $H$ descends to an integrable complex
structure on $N$. It turns $N$ into a compact complex surface. In
this realm, $N$ is known as a Kodaira surface. One of the key
results in \cite{Poon} is the following.
\begin{prop}\label{poon} On the Kodaira surface $N$, the complex structure $J$ and the symplectic structures
\begin{equation}\label{uv symplectic}
u_1(e^1\wedge e^3-e^2\wedge e^4)+v_1(e^1\wedge e^4+e^2\wedge e^3), \quad u_1^2+v_1^2\neq 0,
\end{equation}
are in the same deformation family of generalized complex structures.
\end{prop}
The deformation family could be given explicitly  in terms of a
choice of $(-i)$-eigenspace of an invariant generalized complex
structure. Choose an ordered basis for $(\lie h\oplus \lie
h^*)_\mathbb{C}$ as follows.
\begin{eqnarray*} &\half(e_1+ie_2), \quad
\half(e_3+ie_4), \quad e^1+ie^2, \quad e^3+ie^4,&\\
&\hspace{1.0in}\half(e_1-ie_2), \quad \half(e_3-ie_4), \quad
e^1-ie^2, \quad e^3-ie^4.&
\end{eqnarray*}
Then the $(-i)$-eigenspace is spanned by the row vectors:
\begin{equation}\label{Gamma}
\left(
\begin{array}
[c]{cccccccc}
1 & 0 & 0 & 0 & t_{3} & 0 & 0 & t_{1}\\
0 & 1 & 0 & 0 & 0 & t_{2} & -t_{1} & 0\\
0 & 0 & 1 & 0 & 0 & t_{4} & -t_{3} & 0\\
0 & 0 & 0 & 1 & -t_{4} & 0 & 0 & -t_{2}
\end{array}
\right),
\end{equation}
where $t_1, \dots, t_4$ are complex numbers. When all of them are equal to zero, the distribution is due to the classical
complex structure $J$. When $t_1=t_4=0$, this distribution is due to a generic
classical complex structure. On the other hand, the generalized complex structures determined by
the symplectic form given by (\ref{uv symplectic}) is contained in this family with $t_2=t_3=0$ and
\[
t_1=\frac{i}2(u_1+iv_1), \quad t_4=\frac{2i}{u_1-iv_1}=\frac{1}{{\overline t}_1}.
\]

Note that not all symplectic forms on the Kodaira manifold is contained in the family (\ref{Gamma}). However due to a combination of
(\ref{symplectic family}) with (\ref{Gamma}), the complex structure $J$ and all symplectic forms on $N$ are contained in the same
connected component of generalized deformation family.

\subsection{Geometry on a SO(2)-bundle over a Kodaira surface}
Now we apply the general construction in Section \ref{central
extension}
 to the Kodaira manifold $N$.  Choose the symplectic form
\[
\omega=-(e^1\wedge e^3-e^2\wedge e^4).
\]
Let $M$ be the principal $\SO(2)$-bundle on $N$ with characteristic class $-\omega$. It is covered by a five-dimensional
simply-connected
nilpotent group $G$, which is a central extension of $H$. Let $e_5$ be the fundamental vector field of the principal bundle. Let
$e^5$ be a connection 1-form. Then the structure equations on $\lie g$ are
\begin{equation}
[e_1,e_2]=e_3, \quad [e_1,e_3]=-e_5, \quad [e_2, e_4]=e_5.
\end{equation}
The dual structure equations in terms of the Chevalley-Eilenberg differential are
\begin{equation}
de^3=-e^1\wedge e^2, \quad de^5=-\omega=e^1\wedge e^3-e^2\wedge e^4.
\end{equation}
Treating $e^5$ as a contact 1-form on $G$, we construct its associated generalized contact structure
${\cal J}_1=(F, \eta, \pi, \theta, \varphi)$ as given in Section \ref{classical contact}. We have
\[
F=e_5, \quad \eta=e^5, \quad \pi=e_1\wedge e_3-e_2\wedge e_4,  \quad \theta=e^1\wedge e^3-e^2\wedge e^4,
\quad \varphi=0.
\]
On the other hand, due to the construction of Section \ref{central extension},
 the complex structure $J$ on $H$ induces a generalized contact structure ${\cal J}_0$ with
\begin{eqnarray*}
&F=e_5, \quad \eta=e^5, \quad \pi=0,  \quad \theta=0,&\\
&\varphi=e_2\otimes e^1-e_1\otimes e^2+e_4\otimes e^3-e_3\otimes
e^4.&
\end{eqnarray*}
All  invariant objects on $G$ descend to a co-compact quotient $M$.
As a result of the general construction in Section \ref{central
extension} and Proposition \ref{poon}, we have the following
conclusion.
\begin{prop} The generalized contact structure ${\cal J}_1$ on the manifold $M$ determined by the contact 1-form $e^5$ and the
generalized contact structure ${\cal J}_0$ are in the same deformation family of generalized contact structures.
\end{prop}

Finally, note that the generalized contact structure ${\cal J}_0$ is
not strong in the sense that the space of sections of $L^*$ is not
closed with respect to the Courant bracket on the manifold $M$. One
may check it directly through the given structure equations. One may
also observe that $d\eta=-\omega=e^1\wedge e^3-e^2\wedge e^4$. With
respect to the given $\varphi$, it is type $(2,0)+(0,2)$. Therefore,
the obstruction for the integrability of $L^*$ does not vanish. Due
to Lemma \ref{type 11 lemma}, the triple $(F, \eta, \varphi)$ on $M$
is \em not \rm a normal almost contact structure.

\end{document}